  \documentclass[12pt]{amsart}
\usepackage[utf8]{inputenc}
\usepackage[T2A]{fontenc}
\usepackage{enumerate, amssymb, amsmath, amsthm}
\usepackage{amscd}
\usepackage{graphicx}
\usepackage{hyperref}
\hypersetup{
    colorlinks=true,
    linkcolor=blue,
    filecolor=blue,      
    urlcolor=blue,
    citecolor=magenta
}
 
\usepackage[a4paper,hmargin=2.5cm,vmargin=2.5cm]{geometry}
\usepackage{marginnote}
\usepackage[english]{babel}
\graphicspath{ {images/} }
\usepackage{wrapfig}
\usepackage[matrix,arrow,curve]{xy}\usepackage{mathabx}
\usepackage[OT2, T1]{fontenc}
\DeclareSymbolFont{cyrletters}{OT2}{wncyr}{m}{n}
\DeclareMathSymbol{\Sha}{\mathalpha}{cyrletters}{"58}

\usepackage[a4paper,hmargin=2.5cm,vmargin=2.5cm]{geometry}
\usepackage[english]{babel}
\graphicspath{ {images/} }
\usepackage{wrapfig}
\usepackage[matrix,arrow,curve]{xy}\usepackage{mathabx}
\usepackage{tikz}
\usetikzlibrary{intersections}

\begin{document}
\renewcommand{\refname}{Bibliography}
\newtheorem{prop}{Proposition}[section]
\newtheorem{thrm}[prop]{Theorem}
\newtheorem{lemma}[prop]{Lemma}
\newtheorem{cor}[prop]{Corollary}
\newtheorem{mainthm}{Theorem}
\newtheorem{maincor}[mainthm]{Corollary}
\theoremstyle{definition}
\newtheorem{df}[prop]{Definition}
\newtheorem{ex}[prop]{Example}
\newtheorem{rmk}[prop]{Remark}
\newtheorem{conj}[prop]{Conjecture}
\newtheorem{cl}[prop]{Claim}
\newtheorem{q}[prop]{Question}
\newtheorem{constr}{Construction}
\renewcommand{\proofname}{\textnormal{\textbf{Proof:  }}}
\renewcommand{\refname}{Bibliography}
\renewcommand{\themainthm}{\Alph{mainthm}}
\renewcommand{\themaincor}{\Alph{maincor}}

\renewcommand{\phi}{\varphi}
\renewcommand{\epsilon}{\varepsilon}

\renewcommand{\C}{\mathbb C}
\newcommand{\Z}{\mathbb Z}
\newcommand{\Q}{\mathbb Q}
\newcommand{\R}{\mathbb R}
\newcommand{\N}{\mathbb N}
\newcommand{\Fp}{\mathbb{F}_p}
\newcommand{\Fq}{\mathbb{F}_q}

\renewcommand{\O}{\mathcal O}
\newcommand{\g}{\mathfrak g}
\newcommand{\h}{\mathfrak h}
\newcommand{\E}{\mathcal E}
\newcommand{\F}{\mathcal F}
\newcommand{\m}{\mathfrak{m}}

\renewcommand{\i}{\sqrt{-1}}
\renewcommand{\o}{\otimes}
\newcommand{\di}{\partial}
\newcommand{\acts}{\lefttorightarrow}
\newcommand{\dibar}{\overline{\partial}}
\newcommand{\im}{\operatorname{im}}
\renewcommand{\ker}{\operatorname{ker}}
\newcommand{\Hom}{\operatorname{Hom}}
\newcommand{\tr}{\operatorname{tr}}
\newcommand{\codim}{\operatorname{codim}}
\newcommand{\rk}{\operatorname{rk}}
\newcommand{\nilp}{\operatorname{nilp}}
\newcommand{\hdot}{{\:\raisebox{3pt}{\text{\circle*{1.5}}}}}
\newcommand{\Supp}{\operatorname{Supp}}
\newcommand{\Alb}{\operatorname{Alb}}
\newcommand{\alb}{\operatorname{alb}}
\newcommand{\Hilb}{\operatorname{Hilb}}
\newcommand{\Sh}{\operatorname{Sh}}
\newcommand{\sh}{\operatorname{sh}}
\newcommand{\CP}{\mathbb{C}\mathbf{P}}
\newcommand{\Isom}{\operatorname{Isom}}
\newcommand{\Sym}{\operatorname{Sym}}
\newcommand{\Stab}{\operatorname{Stab}}
\newcommand{\Aut}{\operatorname{Aut}}
\newcommand{\sslash}{\mathbin{/\mkern-6mu/}}
\newcommand{\Pic}{\operatorname{Pic}}
\newcommand{\V}{\mathbb{V}}

\newcommand{\GL}{\operatorname{GL}}
\newcommand{\SL}{\operatorname{SL}}
\newcommand{\SU}{\operatorname{SU}}
\renewcommand{\U}{\operatorname{U}}
\newcommand{\SO}{\operatorname{SO}}
\newcommand{\Ogr}{\operatorname{O}}
\newcommand{\Sp}{\operatorname{Sp}}

\newcommand{\gl}{\mathfrak{gl}}
\renewcommand{\sl}{\mathfrak{sl}}
\newcommand{\su}{\mathfrak{su}}
\newcommand{\so}{\mathfrak{so}}
\renewcommand{\sp}{\mathfrak{sp}}
\renewcommand{\g}{\mathfrak{g}}
\renewcommand{\h}{\mathfrak{h}}
\newcommand{\z}{\mathfrak{z}}
\newcommand{\ad}{\operatorname{ad}}
\newcommand{\cd}{\operatorname{cd}}
\newcommand{\an}{\operatorname{an}}
\newcommand{\orb}{\operatorname{orb}}
\newcommand{\Trop}{\operatorname{Trop}}
\newcommand{\T}{\mathbb{T}}
\newcommand{\Id}{\operatorname{Id}}
\newcommand{\ab}{\operatorname{ab}}

\newcommand{\tX}{\widetilde{X}}
\newcommand{\tM}{\wideilde{M}}
\newcommand{\tF}{\widetilde{F}}
\newcommand{\G}{\mathcal G}
\renewcommand{\L}{\mathcal L}
\renewcommand{\sp}{\mathrm{sp}}

\binoppenalty = 10000
\relpenalty = 10000

\title{An affine surface with $3$-step nilpotent fundamental group}
\author{Sergey Polevoy}
\author{Vasily Rogov}

\begin{abstract}We construct a smooth complex affine surface $Y$ such that $\pi_1(Y)$ is virtually $3$-step nilpotent, but not virtually $2$-step nilpotent. This gives a negative answer to the quasi-projective version of Campana's nilpotency problem. 
\end{abstract}

\maketitle

\tableofcontents

\section{Introduction}

In \cite{Campana95} Campana constructed an example of a compact K\"ahler manifold with $2$-step nilpotent fundamental group and raised the question of which virtually nilpotent groups can be realised in such a way. This is, of course, a special case of Serre's problem on which finitely presented groups are fundamental groups of compact K\"ahler manifolds (\emph{K\"ahler groups}). In the same paper he showed that such a group $\pi_1(X)$ is virtually at most $2$-step nilpotnet if various additional assumptions are satisfied  (see also \cite{CT}). Whether there exists  nilpotent K\"ahler groups of nilpotency higher than $2$ remains open.

In \cite{AC} Aguilar and Campana raised a similar question for fundamental groups of normal quasi-projective varieties. Several recent results \cite{ Rogomin,  Shim2, Shim1, CDHP, Rogcd} show  under various additional assumptions that if the fundamental group of a normal quasi-projective variety is virtually nilpotent, it is virtually at most step $2$ nilpotent, so one would expect the same to hold in general.

In this paper we show that such an expectation is false.

\begin{mainthm}
There exists a smooth complex affine surface $Y$ such that $\pi_1(Y)$ is virtually $3$-step nilpotent, but not virtually $2$-step nilpotent. More precisely, every finite index subgroup of $\pi_1(Y)$ has rational nilpotency class $3$.
\end{mainthm}

The surface $Y$ is obtained as a birational modification at infinity of an elliptic arrangement on an abelian surface. Under \emph{birational modification at infinity} we mean that one performs some blow-ups centered on the boundary divisor are performed and then adds some of the new bounndary components to the surface. In particular, $Y$ is birational to an open subset of an abelian surface. We describe the construction explicitly in subsection \ref{construction subsec}.

Let us outline relations of this result to the previous works:
\begin{itemize}
\item Hirsch length of $\pi_1(X)$ equals $8=\rk L + \dim V$(see Propositions \ref{metabelian prop} and \ref{i2 prop}) which is the smallest possible for such an example according to \cite{Shim2}. We do not know whether the mixed Hodge structure on the Maltsev completion of $\pi_1(Y)$ is the one constructed in \cite[Section 8]{Shim2}, or, more generally, whether the construction of \emph{loc.cit.} can be realised as the mixed Hodge structure on the Maltsev completion of the fundamental group of a normal quasi-projective variety.
\item The surface $Y$ is not aspherical. Indeed, suppose it is. Pass to a finite \'etale cover $Y' \to Y$ such that $\pi_1(Y')$ is torsion-free $3$ step nilpotent. Surface $Y'$ is still aspherical and affine. By the Andreotti-Frankel theorem it has the homotopy type of an aspherical $2$-dimensional CW complex. But the cohomological dimension of $\pi_1(Y')$ equals its Hirsch length, i.e. $8$. This aligns with the results of \cite{Rogcd}.
\item The surface $Y$ is not special in the sense of Campana. In fact, a computation shows that it is of log general type. This aligns with the results of \cite{CDHP}.
\item Further computations show that $\pi_1(Y)$ is not torsion-free. We not include these computations in this version of the paper.
\item When this project was ready, Taito Shimoji told us that he obtained (using ChatGPT) different examples of quasi-projective varieties with virtually nilpotent fundamental group of rational nilpotency $3$, and even higher. His examples are different from ours, although quite close in spirit, and the Hirsch length is much bigger.
\end{itemize}

The paper is organised in the following way. In section \ref{prelim sec} we collect general facts from geometry, topology and group theory.  In subsection \ref{group prep} we recall the relation between metabelian and nilpotent groups, the specifics in the context of fundamental groups of algebraic varieties, and give an explicit recipe for computing the nilpotency of a metabelian virtually nilpotent group. 

In subsection \ref{exact prep} we give a criterion for $\pi_1$-exactness for a relatively proper fibration of complex manifolds, which is a variation of Nori's lemma. 

In subsection \ref{surgery prep} we introduce a topological lemma that helps to control fundamental groups under birational modifications at infinity and is essentially an application of van Kampen's theorem.

Section \ref{main sec} is the main one. We give an explicit construction of $Y$ in \ref{construction subsec}. We also describe $Y$ as a partial compactification of the special fiber of a certain explicit family of quasi-projective surfaces. This point of view, and the comparison with the fundamental group of the general fiber, play a key role in the proof.

In subsection \ref{local subsec} we study the local behaviour of fundamental groups under birational transformations centered at a point where several boundary components meet.

In subsection \ref{metabelian subsec} we use the local computations of subsection \ref{local subsec} together with the glueing lemma of subsection \ref{surgery prep} to get a presentation of $\pi_1(Y)$. In the course of the proof, we show that it is metabelian, which immediately implies virtual nilpotency by \cite{rogbns} or \cite{AN}.

In subsection \ref{subsec last} we carry out the program of subsection \ref{group prep} and  compute the rational lower central series on $\pi_1(Y)$, proving that it has rational nilpotency class $3$ (Theorem \ref{final thrm}).
\\

\textbf{Acknowledgements.} We are thankful to Rodolfo Aguilar, Philippe Eyssidieux and Taito Shimoji for their interest in this work.  
\\

\textbf{AI Usage Disclosure.} The AI input is substantial. The example was suggested by ChatGPT Astra, as a result of a series of prompts made by the first-named author. After that, the second-named author read and verified the generated construction, simplified the proof, and wrote an exposition that meets his personal criteria for an understandable mathematical text. Additionally, Claude Opus 5.5 was used during the preparation of text for searching references, preparing TikZ illustration, checking for typos and polyshing the style. Finally,the manuiscript was independently verified for flaws using ChatGPT Astra and Claude Opus. The authors bear full responsibility for the correctness of the statements.

\section{Technical preliminaries}\label{prelim sec}

Throughout this paper, for an abelian  group $L$ we denote its character torus $\Hom(L, \C^{\times})$ by $\mathbb{T}(L)$. For a group $G$ we denote by $DG$ the derived subgroup.

\subsection{Preliminaries from group theory}\label{group prep}

Recall that a  group $\Gamma$ is called \emph{metabelian} if it can be written as
\[
1 \to M \to \Gamma \to N \to 1
\]
with $M$ and $N$ abelian. Observe that even when $\Gamma$ is finitely presented, $M$ need not be finitely generated.

\begin{thrm}\label{meta-thrm}
Let $\Gamma$ be a finitely presented metabelian group,
\[
1 \to M \to \Gamma \to L \to 1.
\]
Then $\Gamma$ is solvable. If it is moreover virtually nilpotent, the following hold:
\begin{enumerate}
\item the group $M$ has finite rank;
\item Let $V:=M \o \Q$ considered as an $L$-module and $V_{\C}=V \o \C$. Decompose its semisimplification into the sum of characters
\[
V_{\C}^{ss}=\bigoplus_i V_{\chi_i}, \ \chi_i \in \mathbb{T}(L).
\]
Then each $\chi_i$ is torsion;
\item There exists a finite index subgroup $L_1 \subseteq L$ such that the action of $L_1$ on $V$ is purely unipotent, i.e. $V_{\C}^{ss}$ is trivial as an $L_1$-module. 
\end{enumerate}
\end{thrm}
\begin{proof}
It is clear that $[\Gamma, \Gamma] \subseteq M$, thus $\Gamma$ is solvable. A subgroup of a finitely generated virtually nilpotent group is finitely generated.  For the statement about characters see, e.g., \cite[Lemma 4.4.]{AN}, see also \cite[Proposition 12.12]{Py}. Finally, for $L_1$ one  just takes $\bigcap_i \ker \chi_i$.
\end{proof}
\begin{rmk}
A metabelian nilpotent group need not be $2$-step nilpotent. For example, let $A \in \SL_r(\Z)$ be a unipotent matrix with a single maximal Jordan block. Then $\Z^r \rtimes_{A} \Z$ is a metabelian $r$-step nilpotent group.
\end{rmk}

\begin{prop}\label{triv supp prop}
In the setting of Theorem \ref{meta-thrm} the following conditions are equivalent:
\begin{enumerate}
\item\label{one} $L_1=L$, i.e. the action of $L$ is (purely) unipotent;
\item \label{two} $\operatorname{Supp}(V_{\C}^{ss})=\{1\}$, where $\operatorname{Supp}(V_{\C}^{ss}):=\{\chi \ | V_{\chi} \neq 0 \} \subseteq \mathbb{T}(L)$;
\item \label{three} $I^kV=0$ for some power $k \ge 1$ of the augmentation ideal $I \subseteq \Z[L]$.
\end{enumerate}
\end{prop}
\begin{proof}
$L_1=\bigcap_{\chi_i \in \operatorname{Supp}(V^{ss}_{\C})} \ker \chi_i$, therefore (\ref{one}) is equivalent to (\ref{two}).

Assume that $L_1=L$. Denote by $\rho(\ell)$ the action of $\ell \in L$ on $V$. Since all elements $\rho(\ell), \ \ell \in L$ mutually commute, they have a common invariant flag $0=V_0 \subset \ldots \subset V_n =V_{\C}$. The augmentation ideal is linearly generated by elements of the form $(\ell-1)$ and triviality of the support implies that $(\ell-1)\cdot V_j\subseteq V_{j-1}$ and thus at least $I^nV=0$. On the other hand, if $I^kV=0$ for some $k$, then $(\rho(\ell)-1)^k=0$ and each $\rho(\ell)$ acts unipotently.
\end{proof}

\begin{prop}\label{metab kill prop}
Let $X$ be a smooth complex algebraic variety. Suppose that $\pi_1(X)$ is metabelian. Then it is virtually nilpotent.
\end{prop}
\begin{proof}
This follows from \cite{AN} in the case where the normal abelian subgroup is finitely generated, and from \cite{rogbns} in the general case. 
\end{proof}

For a finitely generated group $\Gamma$ we denote by $\gamma_s\Gamma$ its lower central series:
 \[
\gamma_0\Gamma=\Gamma, \qquad \gamma_s\Gamma=[\Gamma, \Gamma_{s-1}\Gamma].
\]
 Assume that $\Gamma$ is virtually nilpotent and $N \subseteq \Gamma$ a finite index nilpotent subgroup. Define the  \emph{rational nilpotency class of $\Gamma$} as
\[
\nilp_{\Q}(\Gamma)=\nilp(N/\operatorname{Tors}(N))=\min\{s \ | \ \gamma_{s}N \text{ is finite} \}.
\]

This does not depend on the choice of $N$ and is stable under passing to finite index subgroups. In particular, $\Gamma$ is not virtually $s$-step nilpotent for $s<\nilp_{\Q}(\Gamma)$.

Suppose now that in the notation of Theorem \ref{meta-thrm} $\Gamma$ is virtually nilpotent. Suppose also for simplicity that $L_1=L$, i.e. all the characters $\chi_i$ are trivial. We may as well assume that $L$ is torsion free.
Below we give a tool for computing $\nilp_{\Q}(\Gamma)$ in this setting.

The group $M$ carries a structure of a $\Z[L]$-module, induced from the $L$-action by conjugations. Let $I \subset \Z[L]$ be the augmentation ideal. Observe that  $IM$ is a normal subgroup of $\Gamma$ and $M/IM$ is central in $\Gamma/IM$.

Consider the map 
\[
\eta \colon \Lambda^2L \to M/IM, \ \ell_1 \wedge \ell_2 \mapsto [\widetilde{\ell_1}, \widetilde{\ell_2}] \ \operatorname{mod} IM,
\]
where $\widetilde{\ell_i}$ are lifts of $\ell_i$ to $\Gamma/IM$. This is well-defined because two different lifts differ by an element from $M/IM$ that is central in $\Gamma/IM$.

Denote by $P_0$ the image of $\eta \colon \Lambda^2L \to M/IM \subseteq \Gamma/IM$ and by $P$ its preimage in $\Gamma$.
\begin{prop}\label{integral criterion}
In the notation above, the following hold:
\begin{enumerate}
\item \label{first} $\gamma_1\Gamma:=[\Gamma, \Gamma]$ is generated by $P$;
\item \label{next} $\gamma_k\Gamma=I\cdot\gamma_{k-1}\Gamma$ for $k\ge 2$ (this makes sense because $\gamma_k\Gamma \subseteq M$).
\end{enumerate}
\end{prop}
\begin{proof}
(\ref{first}) Take any element of the form $[g_1, g_2]$ in $\Gamma$. Let $\pi \colon \Gamma \to \Gamma/IM$ be the projection. Then $\pi([g_1, g_2])=[\pi(g_1), \pi(g_2)]$ is a commutator in $\Gamma/IM$. But since $M/IM$ is central in $\Gamma/IM$, every commutator in $\Gamma/IM$ is contained in $P_0$. Therefore, $\pi([g_1, g_2]) \in P_0$ and $[g_1, g_2] \in P$.
\\

(\ref{next}). By the previous item, $\gamma_1\Gamma \subseteq M$, thus $\gamma_{k-1}\Gamma \subseteq M$ for $k \ge 2$.  Let $x \in \gamma_{k-1}\Gamma$ and $g \in \Gamma$. The commutator $[g,x]$ depends only on the image of $g$ in $L$.  On the other hand, for $\ell \in L$ with a lift $\widetilde{\ell}$ and any  $ x \in M$ holds $[\widetilde{\ell}, x]=a_{(\ell-1)}(x)$, where  $a_r \colon M \to M$ denotes the action of an element $r \in \Z[L]$. The augmentation ideal $I$ is generated by elements of the form $(\ell-1), \ \ell \in L$. We conclude that $\gamma_k\Gamma$ is generated by $I \cdot \gamma_{k-1}\Gamma$. The latter is a subgroup, so the claim follows.
\end{proof}

\begin{lemma}\label{rational criterion}
Let 
\[
1 \to M \to \Gamma \to L \to 1
\]
be a metabelian virtually nilpotent group. Let $V:=M \otimes \Q$. Let $I \subset \Z[L]$ be the augmentation ideal and assume that one of the conditions from Proposition \ref{triv supp prop} holds.  As before, denote by $\pi \colon \Gamma \to \Gamma/IM$ the projection, by $P_0$ the image of the commutator map $\Lambda^2L \to M/IM$ and  $P=\pi^{-1}(P_0)$.  The action of $\Z[L]$ on $M$ extends to $V$ by linearity. Let $W$ be the rational span of $P$ in $V$. For a subspace $V_1 \subseteq V$ denote by $\nu(V_1)$ the minimal power of $I$ annihilating $V_1$. Then the following hold:
\begin{enumerate}
\item \label{clear} $IV \subseteq W \subseteq V$;
\item \label{almost clear} $\nilp_{\Q}(\Gamma)=\nu(W)+1$;
\item \label{pm1} $\nu(V) \le \nilp_{\Q}(\Gamma) \le \nu(V)+1$ and $\nilp_{\Q}(\Gamma)=\nu(V)$ if $V=W$.
\end{enumerate}
\end{lemma}
\begin{proof}
(\ref{clear}) follows from the fact that $IV$ is generated by $IM$ and $IM$ is contained in $P$. From Proposition \ref{integral criterion} it follows that the group $\gamma_k\Gamma$ rationally spans $I^{k-1}W$. This implies (\ref{almost clear}). The last item (\ref{pm1}) is immediate from the first two.
\end{proof}
\subsection{A $\pi_1$-exactness lemma}\label{exact prep}

In this subsection we provide a criterion for $\pi_1$-exactness of a holomorphic fibration on a complex manifold.

\begin{lemma}\label{exactness lemma}
Let $X$ be a smooth complex manifold with $\dim X \ge 3$. Let $D = D_1 \cup \ldots \cup D_m$ be a simple normal crossing divisor on $X$. Let $f \colon X \to C$ be a  proper holomorphic submersion with connected fibers to a complex curve $C$ of genus\footnote{We do not require $C$ to be compact, but we only need to exclude the case $C=\CP^1$.} at least $1$.  Let $X^{\circ}=X \setminus D$ and denote by $F \subset X$ a general fiber and $F^{\circ}=F \cap X^{\circ}$.

Suppose that the following hold:
\begin{itemize}
\item For any positive dimensional stratum $D_J=D_{i_1} \cap \ldots \cap D_{i_k}$ the map $f|_{D_J}$ is a submersion;
\item for a zero dimensional stratum $\{p\}=D_{i_1} \cap \ldots \cap D_{i_n}$ there exist local coordinates centered at $p$ such that $D$ is locally given by $\{z_1 \cdot \ldots \cdot z_n=0\}$ and the function $f$ takes the form
\[
f(z)=f(p)+\sum_{i} a_iz_i, \ a_i\neq 0  \ \forall i.
\]
\end{itemize}
Then there is an exact sequence
\[
1 \to \pi_1(F^{\circ}) \to \pi_1(X^{\circ}) \to \pi_1(C) \to 1.
\]
\end{lemma}
\begin{proof}
Exactness at the right and the middle terms  follows from Nori's lemma \cite[Lemma 1.5]{Nori}. What is left to check is the injectivity of $\pi_1(F^{\circ}) \to \pi_1(X^{\circ})$.

For each finite dimensional stratum $\{p_i\}=\bigcap_{j \in J_i} D_j$ choose a small ball $B_i \subseteq X$ around $p_i$ that does not contain any other $p_j$. Set $q_i=f(p_i)$ and let $\Delta_i\subset C$ be a small disk containing $f(B_i)$. Let $E_i:=f^{-1}(\Delta_i)$.

Choose a fiber $F_t=f^{-1}(t)$ that does not contain any of $p_j$'s. Suppose first that $t \notin \bigcup_i \Delta_i$.

By Thom's first isotopy lemma  in the form of \cite[Proposition 11.1]{Mather}, the restriction of $f$ gives a locally trivial fibration $X \setminus (\bigcup_i E_i) \to C \setminus (\bigcup_i \Delta_i)$ . Since $\pi_2(C\setminus (\bigcup_i \Delta_i))=0$, the map $\pi_1(F_t^{\circ}) \to \pi_1(X^{\circ} \setminus \bigcup_i E_i)$ is injective.

Now let us deal with the case $t \in \Delta_i$. Let $Z_t:=F_t \cap B_i$ be the local nearby fiber. From the explicit local form $f(z)=f(p_i) + \sum_i a_iz_i$ we know that $B_i \cap X^{\circ}$ is homeomorphic to the complement of the coordinate hyperplanes in the unit ball in $\C^n$, i.e. has homotopy type of an $n$-dimensional torus $(S^1)^{\times n}$. At the same time, $Z_t\cap X^{\circ}$ is homeomorphic to the complement of $n$ hyperplanes in general position in $\C^{n-1}$ intersected with a ball. The homotopy type of the pair $(B_i \cap X^{\circ}, Z_t \cap X^{\circ})$ is described by Hattori theorem \cite{Hat} that says that it is $(n-1)$-connected. Thus $\pi_1(Z_t \cap X^{\circ}) \to \pi_1(B_i \cap X^{\circ})$ is injective.

The statement now follows from these two arguments combined with repeated application of  van Kampen's theorem.
\end{proof}

\begin{cor}\label{2-exact cor}
Let $(X, D)$ be a simple normal crossing pair, $\dim X \ge 4$. Let $f_1 \colon X \to C_1$ and $f_2 \colon X \to C_2$ be holomorphic submersions onto irrational complex curves. Denote by $f=(f_1, f_2) \colon X \to C_1 \times C_2$ the map to the product. For a point $(t_1, t_2) \in C_1 \times C_2$ denote  $X_{(t_1, t_2)}:=f^{-1}(t_1,t_2)$ and $D_{t_1, t_2}:=D \cap X_{t_1, t_2}$. Assume that the following hold:
\begin{itemize}
\item  $f_1$ satisfies the assumptions of Lemma \ref{exactness lemma} on $(X, D)$;
\item for a generic fiber $X_{t_1}:=f_1^{-1}(t_1)$ the map $f_2|_{X_{t_1}} \colon X_{t_1} \to C_2$ satisfies the assumptions of Lemma \ref{exactness lemma} on the pair $(X_{t_1}, D_{t_1}:=X_{t_1} \cap D)$.
\end{itemize}
Then for a generic point $(t_1, t_2) $ there is an exact sequence
\[
1 \to \pi_1(X_{t_1, t_2} \setminus D_{t_1, t_2}) \to \pi_1(X \setminus D) \to \pi_1(C_1 \times C_2) \to 1.
\]
\end{cor}
\begin{proof}
The surjectivity of $\pi_1(X \setminus D) \to \pi_1(C_1 \times C_2)$ follows from the fact that $(f_1)_*$ is surjective  and its kernel surjects onto $\pi_1(C_2)$ by diagram chasing. 

Let us check that $\pi_1(X_{t_1, t_2} \setminus D_{t_1, t_2})$ surjects on the kernel of $f_*$. Suppose that $\gamma \in \ker f_*$. Then $\gamma \in \ker (f_1)_*$ and from the exact sequence
\[
1 \to \pi_1(X_{t_1} \setminus D_{t_1}) \to \pi_1(X \setminus D) \to \pi_1(C_1)\to 1
\]
there exists a lift $\widetilde{\gamma}$ of $\gamma$ to $\pi_1(X_{t_1} \setminus D_{t_1})$. 

At the same time, there is an exact sequence
\[
1 \to \pi_1(X_{t_1, t_2} \setminus D_{t_1,t_2}) \to \pi_1(X_{t_1} \setminus D_{t_1}) \to \pi_1(C_2) \to 1 .
\]
The element $\widetilde{\gamma}$ lies in the kernel of the map to $\pi_1(C_2)$ because this map factors through $\pi_1(X\setminus D)$ and $\gamma \in \ker(f_2)_*$, therefore it comes from $\pi_1(X_{t_1, t_2} \setminus D_{t_1, t_2})$.

Finally, the map $\pi_1(X_{t_1, t_2} \setminus D_{t_1,t_2}) \to \pi_1(X \setminus D)$ is injective as a composition of injective maps $\pi_1(X_{t_1, t_2} \setminus D_{t_1,t_2}) \to \pi_1(X_t \setminus D_t)$ and $\pi_1(X_t \setminus D_t) \to \pi_1(X \setminus D)$.
\end{proof}

\subsection{A local modification lemma}\label{surgery prep}

In this subsection we prove a lemma that describes  the behaviour of $\pi_1$ under local modifications that is an easy application of van Kampen's theorem. First, we need a group-theoretic proposition.

\begin{prop}\label{vankampen prop}
Assume we are given a commutative diagram of groups
\[
\xymatrix{
C \ar[r]^{i_A} \ar[d]_{i_B} \ar[rrd]^{\gamma} & A \ar[rrd]^{\alpha} & & \\
B \ar[rrd]_{\beta} & & C' \ar[r]_{i_{A'}} \ar[d]_{i_{B'}} & A' \\
& & B' &
}
\]
such that  $\alpha, \beta$ and $\gamma$ are surjective.  Let $G:=A*_CB$ and $G'=A'*_{C'} B'$ . Then  the induced map $\phi \colon G \to G'$ is surjective  and its kernel is normally generated by the images of $\ker \alpha$ and $\ker \beta$.
\end{prop}
\begin{proof}
First, $\phi$ is surjective because $G'$ is generated by $A'=\im \alpha = \im \phi|_A$ and $B'=\im \beta = \im  \phi|_B$.

Denote by $N \trianglelefteq G$ the normal subgroup of $G$ generated $\ker\alpha$ and $\ker \beta$. Clearly, $N$ is contained in $\ker \phi$, thus $\phi$ descends to a surjective map $\overline{\phi} \colon G/N \to G'$. Each $a' \in A'$ can be lifted to an element $a \in A$ and the image of $a$ under $A \to G \to G/N$ does not depend on the lift. Similarly, for $b' \in B'$. This gives two maps $A' \to G/N, \ B' \to G/N$ that agree on $i_{A'}(C')$ and $i_{B'}(C')$ (this follows from the  fact that $G \to G/N$ agrees on $i_A(C)$ and $i_B(C)$ and $\gamma$ is onto). By the universal property of amalgamated product this defines a map $G' \to G/N$ which is an inverse of $\overline{\phi}$. Hence $\overline{\phi} \colon G/N \to G'$ is an isomorphism.
\end{proof}

This proposition has the following topological corollary. Informally speaking, it describes the fundamental group changes under local modifications at  two points, provided that these modifications are surjective on local fundamental groups (e.g. partial compactifications); we do not make the word ``'modification'' precise here.

\begin{lemma}\label{triple cover lemma}
Let $f \colon X \to \widehat{X}$ be a continuous map of CW-complexes. Assume that there is an open cover by connected sets $X=V_1 \cup W \cup V_2$ with pairwise connected intersections $J_i=V_i \cap W, \ i=1,2$ and $V_1 \cap V_2=\varnothing$, and a similar cover $\widehat{X}=\widehat{V}_1 \cup \widehat{W} \cup \widehat{V}_2$ with $\widehat{J}_i=\widehat{V}_i \cap \widehat{W}$. Assume that $f$ maps $V_i, \ W, \ J_i$ into $\widehat{V}_i, \ \widehat{W}, \ \widehat{J}_i$ respectively, and that  moreover the following hold:
\begin{itemize}
\item $J_i \to \widehat{J}_i$ for $i=1,2$ and $W \to \widehat{W}$ are homotopy equivalences;
\item $f_* \colon \pi_1(V_i) \to \pi_1(\widehat{V}_i)$ are surjective for $i=1,2$.
\end{itemize}

Then the map $f_* \colon \pi_1(X) \to \pi_1(\widehat{X})$ is surjective and its kernel is normally generated by the images of $\ker[\pi_1(V_i) \to \pi_1(\widehat{V}_i)]$ for $i=1,2$.
\end{lemma}
\begin{proof}
Define $W_1:=W \cup V_2$ and $\widehat{W}_1:=\widehat{W} \cup \widehat{V}_2$. By van Kampen's theorem, 
\[
\pi_1(X)=\pi_1(W_1) *_{\pi_1(J_1)} \pi_1(V_1)conv
\]
and similarly 
\[\pi_1(\widehat{X})=\pi_1(\widehat{W}_1)*_{\pi_1(\widehat{J}_1)} \pi_1(\widehat{V}_1)
\]
We apply Proposition \ref{vankampen prop} and conclude that the kernel is normally generated by images of kernels of  $(f|_{V_1})_*$ and $(f|_{W_1})_*$. What is left is to determine $\ker[\pi_1(W_1) \to \pi_1(\widehat{W}_1)]$. For this we apply van Kampen's theorem and Proposition \ref{vankampen prop} once again, now for the open covers $W_1=W \cup V_2$ and $\widehat{W}_1=\widehat{W} \cup \widehat{V}_2$. We conclude that the kernel of the map
\[
\pi_1(W_1)=\pi_1(W)*_{\pi_1(J_2)} \pi_1(V_2) \to \pi_1(\widehat{W})*_{\pi_1(\widehat{J}_2)}\pi_1(\widehat{V}_2)=\pi_1(\widehat{W}_1)
\]
is normally generated by the images of $\ker[\pi_1(V_2) \to \pi_1(\widehat{V}_2)]$ and $\ker[\pi_1(W) \to \pi_1(\widehat{W})]$. The latter is trivial.
\end{proof}
\section{The construction}\label{main sec}

\subsection{The surface}\label{construction subsec}
We present the construction of the surface. The notation introduced in this subsection is used until the end of the paper.

The general construction is the following. We start with an elliptic curve $E$, an abelian surface $A=E \times E$ and a divisor $B \subset A$. Then we replace $A$ with a birational modification $\sigma \colon X \to A $ that is an isomorphism outside $B$. We set $\widetilde{B}$ to be the full preimage of $B$ in $X$ and take a divisor $D \subseteq X$ to be a union of some (not all) components of $\widetilde{B}$. The required surface is obtained as $Y:=X \setminus D$.

Now we present the details. Let $E$ be an elliptic curve and $\tau \in E[3]$ a non-zero $3$-torsion element. Let $A:=E \times E$. We consider the following $4$ curves in $A$:
\begin{itemize}
\item $B_1:=\{(2x, x) \ | x \in E\}$;
\item $B_2:=\{(x, 2x) \ | x \in E\}$;
\item $B_{0}:=E \times \{0\}$;
\item $B_{\tau}:=\{\tau\} \times E$.
\end{itemize}
These are four elliptic curves on $A$ that meet in exactly two triple points, $p=(0,0)$ and $q=(\tau, 2\tau)$. Denote $B=B_0 \cup B_1 \cup B_2 \cup B_{\tau}$ and let $A':=\operatorname{Bl}_{p,q}A$ be the blow-up of $A$ in the points $p$ and $q$. Let $\widetilde{B}$ be the full preimage of $B$ in $A'$, i.e. 
\[
\widetilde{B}=B_0' \cup B_1' \cup B_2' \cup B_{\tau}' \cup E_p \cup E_q
\]
where $B'_j$ are strict transforms of $B_j$ and $E_p, \ E_q$ are the exceptional divisors.

Let $p':=B'_1 \cap E_p$ and $q':=B'_2 \cap E_q$. Let $X:=\operatorname{Bl}_{p', q'}A'$ be the  blow-up of $A'$ centered in $p'$ and $q'$ . Let $D_0, \ D_1, \ D_2$ and $D_{\tau}$ be the strict transforms of $B'_0, \ B'_1, \ B'_2$ and $B'_{\tau}$ and $E'_p, \ E'_q$ the strict transforms of $E_p$ and $E_q$.

Let 
\[
D=D_0 \cup D_1 \cup D_2 \cup D_{\tau} \cup E'_p \cup E'_q.
\]
The surface we are interested in is 
\[
Y:=X \setminus D.
\]
\begin{figure}
  \centering
  \includegraphics[width=0.7\textwidth]{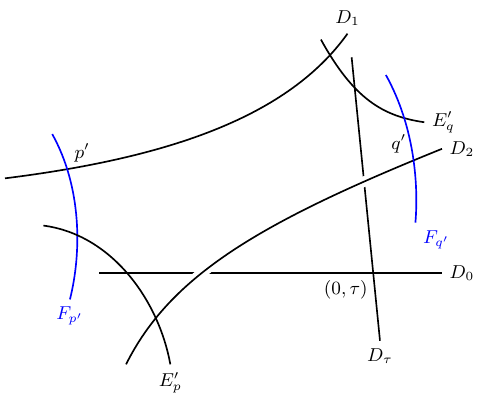}

  \caption{The configuration of components of $D$. The curves $F'_p$ and $F'_q$ are exceptional divisors of the blow-up $X \to A'$; they belong to $Y=X \setminus D$ and each one meets $D$ in two points.}
  \end{figure}
\begin{prop}\label{affine}
The surface $Y$ is affine and of log-general type.
\end{prop}
\begin{proof}
Since $X$ is projective, it is enough to check that $D$ supports an ample divisor.

To apply Nakai-Moishezon criterion,  we first need to verify that $D$ meets every curve in $X$, or, equivalently, there is no complete curve in $Y$. Let $C \subset X$ be a curve that does not intersect $D$. Denote by $\sigma \colon X \to A$ the blow-up map and let $C_0=\sigma(C)$. Since the exceptional locus of $\sigma$ is contained in $D$, the curve $C$ is not contracted by $\sigma$. Observe that $B$ is ample on $A$, so $C_0$ necessary intersects $B$. It can only meet it at points $p$ and $q$ (otherwise, $C$ intersects $D$), and it pass through both of these points (otherwise either $C_0 \cdot B_0$ or $C_0\cdot B_{\tau}$ vanishes, which means that $C_0$ is of the form $E \times \{t_0\}$ or $\{t_0\} \times E$ which creates new intersections of $C_0$ with $B$). Let $\mu_p$ and $\mu_q$ be the multiplicites  of $C_0$  at the points $p$ and $q$ and $\nu_p$ and $\nu_q$ be the multiplicites of its strict transform at  $p'$ and $q'$. Then
\begin{equation}\label{sigma}
\sigma^*C_0=C+\mu_p(E'_p+F_{p'})+\nu_pF_{p'}+\mu_q(E'_q+F_{q'})+\nu_qF_{q'}.
\end{equation}
Let us make use of the fact that $C \cdot D_j=0$ for each component $D_j \subset D$. The total transforms $E'_p+F_{p'}, \ E'_q+F_{q'}$ and the curves $F_{p'}, \ F_{q'}$ are orthogonal to $\sigma^*\operatorname{NS}(A)$ and to each other, and each has square $-1$. In particular $(E'_p)^2=(E'_q)^2=-2$ and $E'_p\cdot F_{p'}=E'_q \cdot F_{q'}=1$. Intersecting (\ref{sigma}) with $E'_p$ gives $0=\mu_p-\nu_p$, so $\mu_p=\nu_p$. For similar reasons, $\mu_q=\nu_q$. Further, 
\[
0=C \cdot D_1=\sigma^*C_0 \cdot \sigma^*B_1-\mu_p-\mu_q - \nu_p,
\]
thus $C_0 \cdot B_1=2\mu_p+\mu_q$. Computations in the same spirit show $C_0 \cdot B_0=\mu_p, \ C_0 \cdot B_{\tau}=\mu_q$ and $C_0 \cdot B_2=\mu_p+2\mu_q$.

Let $\Delta\subset A=E\times E$ be  the diagonal. In the N\'eron-Severi group $NS(A)$ classes of $[B_1]$ and $[B_2]$ decompose as $[B_1]=2[B_0]-[B_{\tau}]+2[\Delta]$ and $[B_2]=2[B_{\tau}]-[B_0]+2[\Delta]$.  This implies that $[B_1]-[B_2]=3([B_0]-[B_{\tau}])$ and
\[
C_0 \cdot (B_1-B_2)=(\mu_p-\mu_q)=C_0 \cdot(3B_0-3B_1)=3(\mu_p-\mu_q).
\]
Therefore, $\mu_p=\mu_q=:\mu$.

Since $K_X=E'_p+2F_{p'}+E'_q+2F_{q'}$ (the sum of the total transforms of the exceptional curves) and $\nu_p=\mu_p=\nu_q=\mu_q=\mu$, we can write $C=\sigma^*C_0-\mu K_X$. Remark that $K_X^2=-4$ and $K_X$ is orthogonal to $\sigma^*C_0$, so $C^2=C_0^2-4\mu^2$ and $K_X\cdot C=4\mu$. By adjunction,
\[
2p_a(C)-2=C\cdot(K_X+C)=C_0^2-4\mu(\mu-1),
\]
therefore
\[
C_0^2=2p_a(C)-2+4\mu(\mu-1)\geq 4\mu(\mu-1),
\]
since $C$ is a curve on an abelian surface, so it is not rational and $p_a(C)\geq1$.

Antoher inequality on $C_0^2$ is given by Hodge index theorem on $A$. Consider the divisor $H:=B_0+B_{\tau}$. It is ample on $A$, so
\[
C_0^2H^2=2C_0^2 \le (C_0\cdot H)^2=4\mu^2,
\]
hence $C_0^2 \le 2\mu^2$.  The equality would force $[C_0]=\mu [H]$ which is excluded by $C_0 \cdot B_1=3\mu$ while $\mu H \cdot B_1=5\mu$.

We get that 
\[
4\mu(\mu-1) < 2\mu^2
\]
which implies $\mu <2$. As we already explained, $\mu$ is positive. We are left with the only possibility $\mu=1$. Then $C_0^2<2\mu^2=2$  and $C_0^2$ has to be even, so $C_0^2=0$ and $C_0$ is a subgroup. On the other hand $C_0$ is tangent  to $B_1$ at $0$, which forces then $C_0=B_1$ and leads to a contradiction.

To apply Nakai-Moishezon criterion on $X$, it remains to check it on curves supported on $D$, i.e. verify that there exists a divisor $G$ supported on $D$ such that $G \cdot Z>0$ for any irreducible component $Z$ of  $D$. It is enough to take $G=2(D_0+D_1+D_2+D_{\tau}) + E'_p+E'_Q$. This concludes the proof of afiness.

To see that $Y$ is of log general type, observe that $2D=G+E'_p+E'_q$ is a sum of an ample divisor and an effective one, hence it is big. Therefore, $D$ is big and $K_X+D$ is a sum of an effective and a big divisor. This implies that $K_X+D$ is big and $\bar\kappa(Y)=2$.
\end{proof}

In the rest of the paper we are going to study the fundamental group of $Y$. Of course, the main tool of doing this is by comparing with the fundamental group of the surface $U:=A \setminus B$, which $Y$ partially compactifies. The surface $U$ itself has complicated topology that comes from the presence of the triple points $p$ and $q$ on the boundary divisor. To deal with this, we realise $U$ as a special fiber $U_b$ of a family of affine surfaces $U_t$  whose boundary divisors have generically only double points. The construction is as follows.

We start from a $4$-fold $\mathcal{A}:=A \times A =E^{\times 4}$ and an open subset $\mathcal{U}:=(E\setminus\{0\})^{\times 4} \subset \mathcal{A}$. Denote by $\mathcal{B}:=(A \times A) \setminus \mathcal{U}$ the complement.  The notation might look slightly awkward, but the idea is to  interpret $\mathcal{A}$ as a family of surfaces $A_t$ parametrised by $A$. Each fiber $A_t$ will be abstractly isomorphic to $A$, but the projection will be not the standard one, so that the topology of the divisor $B_t=\mathcal{B}\cap A_t$ varies in $t$ and at the special point $b=(\tau, 2\tau)$ one gets $B_b=B$.

Consider the map 
\begin{align*}
\overline{f} \colon  \mathcal{A}=E^{\times 4} &\longrightarrow E^{\times 2} \\
(x_1, x_2, x_3, x_4) &\longmapsto (x_1 + x_2 - x_3 - x_4,\; x_1 + x_3 - 2x_4).
\end{align*}

We also denote by $f$ its restriction to $\mathcal{U}$.

For a point $t=(t_1,t_2)$ in $E^{\times 2}$ the fiber $A_{t}:=\overline{f}^{-1}(t_1,t_2)$ is abstractly isomorphic to $A=E \times E$ (the isomorphism is given e.g. via the projection on $(x_3, x_4)$). However, the pair $(A_t, B_t:=\mathcal{B} \cap A_t)$ does vary globally, as well as the affine surface $U_t=A_t\setminus B_t$.

\begin{prop}\label{exact app prop}
In the notations as above, the following holds:
\begin{enumerate}
\item \label{exact}for a generic $t=(t_1,t_2) \in A$ the  sequence
\[1 \to \pi_1(U_{t}) \to \pi_1(\mathcal{U}) \to \pi_1(A) \to 1\]
is exact;
\item \label{special fiber} for  $b=(\tau, 2\tau)$ there is an isomorphism $A_{b} \xrightarrow{\sim} A$ that maps $B_b$ to $B$. In particular, $U_{b}$ is isomorphic to $U=A \setminus B$ and  embeds in $Y$ as a dense open subset.
\end{enumerate}
\end{prop}
\begin{proof}
(\ref{exact})The components $f_1$ and $f_2$ of $f$ satisfy the assumptions of Corollary \ref{2-exact cor}. The map $\overline{f}$ lifts to a linear map $\C^4 \to \C^2$ given by the matrix
\begin{equation}\label{m}
\widetilde{f}=\begin{pmatrix} 1 & 1 &-1&-1\\1&0&1&-2 \end{pmatrix}
\end{equation}
for which the criteria  are immediate to check.
\\

To prove the second assertion (\ref{special fiber}) set $u=x_3+\tau$ and $v=x_4$. Given 
$$\overline{f}(x_1, x_2, x_3, x_4)=(x_1+x_2-x_3-x_4, x_1+x_3-2x_4)= (\tau, 2\tau),$$ we obtain the expressions
\[
x_1=2\tau-(u-\tau)+2v=2v-u
\]
and 
\[
x_2=\tau-(2v-u)+(u-\tau)+v=2u-v.
\]
In other words, the isomorphism $A_{(\tau, 2\tau)} \xrightarrow{\sim} A$ given by $(u,v)$ sends the four components of the boundary divisor \[
B_{(\tau, 2\tau)}=A_{(\tau, 2\tau)} \cap \left ( \bigcup_i \{x_i=0\} \right)
\]
to the curves $\{2v=u\}, \ \{2u=v\}, \ \{u=\tau\}, \ \{v=0\}$, i.e. to $B_1, \ B_2, \ B_{\tau}, \ B_0$ respectively.

\end{proof}

\subsection{Local topology}\label{local subsec}

Now we are going to study what happens topologically  around the points $p$ and $q$ on $A$ (or, to be more precise, around the points $(p,b)=(0,0,\tau,2\tau)$ and $(q,b)=(\tau, 2\tau, \tau, 2\tau)$ in the total space of the family $\mathcal{A}$). 

The local model is the following. Consider the three planes $Q_1, \ Q_2, \ Q_3$ in $\C^3$ given by the equations $\{x=0\}, \ \{y=0\}$ and $\{x+y+z=0\}$ respectively. For each value $c \in \C$ we get three lines $L_{i,c}:=Q_i \cap P_c$ on the plane $P_c=\{z=c\}$ that have only double intersections for $c \neq 0$ and a single triple point at $c=0$.

We denote $V:=\C^3 \setminus (Q_1 \cup Q_2 \cup Q_3)$ and $V_0:=V \cap P_0$. One can think about $V$ as the total space of a family of affine surfaces $V_z=V \cap P_z$ degenerating at $z=0$.

Now let $X_0$ be the blow up of the plane  $P_0=\{z=0\}$ at the origin. Let $E$ be the exceptional divisor and $M_i$ be the strict transforms of the lines $L_{i,0}$ for $i=1, \ 2, \ 3$. Take the point $\{r\}:=E \cap M_1$ and blow up $X_0$ at this point. We get a new surface $X_1$ with an exceptional divisor $F$ and strict transforms $E', \ M'_i$ of $E$ and $M_i$.

Let $Y_0:=X_1 \setminus \left(E' \cup M'_1 \cup M'_2 \cup M'_3 \right)$.

The reader most probably recognises the set-up: $V$ is  a local model for the intersection of the neighbourhood of $(p,b)$ in $\mathcal{A}$  with $\mathcal{U}$; the special fiber $V_0$ is its intersection with $U_0$ and $Y_0$ is the local model for $Y$ at this point (the geometry around $(q, b)$ is the same).

Our aim is to compute the fundamental groups $\pi_1(V_0), \ \pi_1(V)$ and $\pi_1(Y_0)$ and the natural maps between them.

\begin{prop}\label{local to nearby prop}
There are  group  isomorphisms $\pi_1(V) \simeq \Z^3$ and $\pi_1(V_0) \simeq \Z \times F_2$, where $F_2$ is the free group on two generators. The natural map $\iota_* \colon \pi_1(V_0) \to \pi_1(V)$ is surjective and coincides with the abelianisation.
\end{prop}
\begin{proof}

The isomorphism $\pi_1(V) \simeq \Z^3$ is clear because of an isomorphism $V \simeq (\C^{\times})^3$.

For the second isomorphism observe that the  surface $V_0$ admits a structure of a principal holomorphic $\C^{\times}$-bundle over $\CP^1 \setminus \{0,1,\infty\}$. Every such bundle is topologically trivial, thus $\pi_1(V_0)=\pi_1(\C^{\times})\times \pi_1(\CP^1 \setminus \{0,1,\infty\}) = \Z \times F_2$.

More explicitly, the group $\pi_1(V_0)$ is generated by the two meridians $\lambda_1$ and $\lambda_2$ around $L_1$ and $L_2$ and a central element $\epsilon$ represented by the loop around the origin in a line $L \subset \C^2, \ L \neq L_{i,0}$. If $e_1, e_2, e_3$ are generators of $\Z^3=\pi_1(V)$ with each $e_i$ being represented by a loop around $Q_i$, the map $\iota_*$ is given by
\[
\lambda_1 \longmapsto e_1, \qquad
\lambda_2 \longmapsto e_2, \qquad
\epsilon \longmapsto e_1 + e_2+e_3.
\]

\end{proof}
\begin{prop}\label{local factor prop}
The group $\pi_1(Y_0)$ is isomorphic to $ \Z^2$ and the map $\pi_1(V_0) \to \pi_1(Y_0)$ is surjective and factors through $\iota_*$. 
\end{prop}
\begin{proof}
First of all, recall two  simple general facts about the topology of complements of simple normal crossing divisors:
\begin{enumerate}
\item Let $X$ be a complex manifold and $D=D_1 \cup \ldots \cup D_k$ be a simple normal crossing divisor. Let $D':=D_1 \cup \ldots \cup D_r$ be the union of some of the irreducible components. Then the map $\pi_1(X \setminus D) \to \pi_1(X \setminus D')$ is surjective and its kernel is normally generated by the classes of  meridian loops around $D_j, \ r \le j \le k$;
\item Let $B$ be a complex ball and $D=D_1 \cup D_2$ a simple normal crossing divisor consisting of two irreducible components intersecting at the origin. Let $\gamma_i$ denote the classes of the meridian loops around $D_i$. Let $\pi \colon \widetilde{B} \to B$ be the blow-up at the origin. Of course this does not change the complement $B \setminus D=\widetilde{B} \setminus \pi^*D$. If $\delta$ is the class of the meridian loop around the exceptional divisor in $\pi_1(\widetilde{B} \setminus \pi^*D)$, it equals
\[
\delta = \gamma_1 \gamma_2
\]
(for a suitable choice of orientations).
\end{enumerate}

Now we apply both facts to the construction of $Y_0$. Let $\lambda_i$ be the classes of meridian loops around $L_i$. Then the loop around the exceptional divisor $E$ in 
\[
\pi_1(V_0)=\pi_1(X_0 \setminus (E \cup M_1 \cup M_2 \cup M_3))
\]
is given by the central element $\epsilon$
(see  the proof of Proposition \ref{local to nearby prop}). Applying the same principle once again, we conclude that the class of the meridian loop around the exceptional divisor $F$ equals 
\[
\zeta:=\lambda_1\epsilon.
\]
Thus, $\pi_1(Y_0)$ is the quotient of $\Z\times F_2$ by the subgroup normally generated by $\zeta=\lambda_1\epsilon$, where $\epsilon$ is the central element and $\lambda_1$ one of the two free generators.

The map $j_{*} \colon \pi_1(V_0) \to \pi_1(Y_0)$ is surjective, so $\pi_1(Y_0)$ is generated by $j_*(\epsilon)$ and $j_*(\lambda_2)$. These elements commute, so $\pi_1(Y_0)$ is commutative. Therefore $\pi_1(V_0) \to \pi_1(Y_0)$ factors through $\pi_1(V_0)^{\ab}$ and $\pi_1(V_0)^{\ab} \to \pi_1(Y_0)$ is given by taking the quotient by the subgroup generated by $\langle2e_1+e_2+e_3\rangle \subseteq \Z\langle e_1, \ e_2, \ e_3 \rangle$.  This vector is primitive and the quotient is free of rank 2.

\end{proof}

\subsection{$\pi_1(Y)$ is metabelian}\label{metabelian subsec}
In this subsection we establish that $\pi_1(Y)$ is metabelian. After that, using the results of subsection \ref{group prep} we will be able to conclude that $\pi_1(Y)$ is virtually nilpotent and compute its rational nilpotency class.

Our next goal is to gain a better understanding of the fundamental group of the general fiber $U_t$ and its relation to $\pi_1(U_0)$, and eventually to $\pi_1(Y)$.

The group $\pi_1(U_t)$ can be described as the kernel of the map $f_* \colon \pi_1(\mathcal{U}) \to \pi_1(A)$. Note that $\pi_1(\mathcal{U})=\pi_1(E\setminus\{0\})^{\times 4}\simeq (F_2)^{\times 4}$, where $F_2$ is the free group on two generators.

The map $f_* \colon \pi_1(U) \to \pi_1(A)$ factors through a homomorphism $\overline{f}_* \colon \pi_1(A \times A) \to \pi_1(A)$. 

The map $\overline{f}_*$ is easy to compute; we state the result of a computation as a separate proposition for future reference. 

\begin{prop}
Let $H=\pi_1(E)$, so that $\pi_1(A)=H^{\oplus 2}$ and $\pi_1(A \times A)=H^{\oplus 4}$. Then $\overline{f}_*$ is given by the matrix
\begin{equation}\label{matrix}
\Phi:=\begin{pmatrix} \Id_H & \Id_H & -\Id_H& -\Id_H \\ \Id_H & 0 & \Id_H & -2\Id_H\end{pmatrix}
\end{equation}
\end{prop}

The group $\Pi:=\pi_1(U_t)$ fits into an exact sequence
\[
1 \to (DF_2)^{\times 4} \to \Pi \to \ \ker \Phi \to 1
 \]
 where $DF_2$ is the derived group of $F_2$ and $H$ is viewed as $F_2^{\operatorname{ab}}$. The map $\Phi$ is a surjective homomorphism from $\Z^{8}$ to $\Z^4$, so its kernel $L$ is isomorphic to $\Z^4$. We denote $(DF_2)^{\times 4}$ by $K$, so that $\Pi$ has a structure of an extension of a free abelian group $L$ by an infinitely generated group $K$.
\begin{prop}\label{lem:translates}
For $t = (t_1,t_2)$ set
$\theta_1(t) := t_2 + \tau$ and $\theta_2(t) := t_2 - t_1 + 2\tau$. The isomorphism $A_t \to A$,
$(x_1,\ldots,x_4) \mapsto (u,v) := (x_3+\tau, x_4)$, maps $B_t$ to
\[
B_0 \cup B_\tau \cup \{2u - v = \theta_1(t)\} \cup \{u - 2v = \theta_2(t)\}.
\]
In particular, $B_t$ is obtained from $B$ by translating $B_1$ and $B_2$,
and $\theta_1(b) = \theta_2(b) = 0$.
\end{prop}

\begin{proof}
On $A_t$ one has $x_2 = t_2 - x_3 + 2x_4$ and
$x_1 = t_1 - t_2 + 2x_3 - x_4$; substitute $x_3 = u - \tau$, $x_4 = v$.
\end{proof}

The functions $\theta_1$ and $\theta_2$ give local coordinates on the base of our fibration in the neighbourhood of $b$ that are well-suited to control singularities of the boundary divisor $B_{\theta}, \ \theta=(\theta_1(t), \theta_2(t))$. Near $p$ the three components meet at one point only when
$\theta_2 = 2\theta_1$, and near $q$ only when $\theta_1 = 2\theta_2$. We will think of $\theta=(\theta_1, \theta_2)$ as coordinates in a neihgbourhood of $b=(0,0)$ defined on a small open ball $\Theta$ centered at $b$.
Denote by 
\[
\Theta^{*}:=\{2\theta_1 \neq \theta_2, \ 2\theta_2 \neq \theta_1\} \subset \Theta.
\]
This is exactly the set of points $\theta \in \Theta$ such that the divisor $B_{\theta}$ over this point has only simple double self-intersections. For every such $\theta$ by Proposition 3.2 the inclusion of the fiber $U_{\theta} \subset \mathcal{U}$ identifies
$\pi_1(U_{\theta})$ with $\Pi$.

We denote by $U_{\Theta}$ and $A_{\Theta}$ the preimages of $\Theta$ in $\mathcal{U}$ and $\mathcal{A}$ correspondingly.

Choose two small balls $V'(p) \subset V(p) \subset U_{\Theta} $ centered at $p$  and  $V'(q) \subset V(q) \subset U_{\Theta}$ centered at $q$. Denote $W:=U_{\Theta} \setminus (V'(p) \cup V'(q))$ and $J_p:= W \cap V(p)$ (resp. $J(q):=W \cap V(q)$. For $\theta \in \Theta$ we denote by $V'_{\theta}(p), V_{\theta}(p), W_{\theta}, J_{\theta}$ etc. the intersection of a corresponding set with the fiber $U_{\theta}$.Recall that in $\theta$-coordinates the origin corresponds to $b$, so the fiber $U_0$ is the one we ar eventually interested in.

Let $G_p$ and $G_q$ denote the images of $\pi_1(V_0(p))$ and $\pi_1(V_0(q))$ in $\pi_1(U_0)$.

\begin{lemma}\label{prop:specialization}
Suppose $\theta \in \Theta^*$. There is a continuous map
$\varphi \colon U_0 \to U_{\theta}$ such that $\varphi_* \colon \pi_1(U_0) \to
\pi_1(U_{\theta}) \simeq \Pi$ is surjective and its kernel is normally
generated by $[G_p,G_p]$ and $[G_q,G_q]$. Moreover, $\varphi$ maps
meridian loops around each $B_j$ to meridian loops around the
corresponding component of $B_{j,\theta}$.
\end{lemma}

\begin{proof}
 By Thom's first isotopy lemma there is
a diffeomorphism $\Psi \colon W_0 \to W_{\theta}$ that preserves the distances
to $p$ and $q$ on the shells and maps loops around divisours to loops around divisors.

The components of $B$ through $p$ and $q$ locally look like lines through these
points in flat coordinates.

In particular, inside balls $V_0(p)$ and $V_0(q)$ there are well defined central dilatations that preserve the divisor $B$. Choose intermediate balls $V_0'(p) \subset V_0''(p) \subset V_0(p)$ and $V_0'(q) \subset V_0''(q) \subset V_0(q)$. Let $r \colon U_0 \to U_0$ be a smooth map that projects the balls $V_0''(p))$ and $V_0''(q)$ around $p$ and $q$ radially onto
their boundary spheres and is the identity elsewhere. Set
$\varphi := \Psi \circ r$.

We check the assumptions of Corollary \ref{triple cover lemma} for $\varphi$ and the covers
$U_0 = V_0(p) \cup W_0 \cup V_0(q)$ and
$U_{\theta} = V_{\theta}(p) \cup W_{\theta} \cup V_{\theta}(q)$. The map $r$ is homotopic to the
identity by a radial homotopy that preserves $W_0$ and $J_0(\cdot)$.
Hence $\varphi$ restricts to homotopy equivalences $W_0 \to W_{\theta}$ and
$J_0(\cdot) \to J_{\theta}(\cdot)$, and it maps $V_0(\cdot)$ into
$V_{\theta}(\cdot)$. Moreover, $J_0(p) \hookrightarrow V_0(p)$
is a homotopy equivalence. On the other hand, by Propositions \ref{local to nearby prop} and \ref{local factor prop} we know that   $\pi_1(V_{\theta}(p)) \simeq \mathbb{Z}^3$ is generated by loops around the boundary components and 
\[
\pi_1(V_0(p)) \simeq \mathbb{Z}\times F_2
\to \pi_1(V_{\theta}(p)) \simeq \mathbb{Z}^3
\]
is the abelianisation; in particular it is surjective, and its
kernel is the derived subgroup. The same holds at $q$, and Corollary
\ref{triple cover lemma} gives the claim.
\end{proof}

 \begin{cor}
Let $\Gamma:=\pi_1(Y)$. The map $\pi_1(U_0) \to \Gamma$ factors through $\Pi$.
\end{cor}
\begin{proof}
Follows from Lemma \ref{prop:specialization} and Proposition \ref{local factor prop}.
\end{proof}

We have constructed a surjection $\pi \colon \Pi \to \Gamma$. The benefit from it is that the group of the general fiber $\Pi$ is much easier to understand rather than the special fiber group $\pi_1(U_0)$.

Recall that the group $\Pi$ has structure
\[
1 \to K \to \Pi \to \Z^4 \to 1,
\]
where $K$ is isomorphic to $(DF_2)^{\times 4}$ (here $F_2$ is the free group on two letters and $DF_2$ is its derived subgroup).

Moreover, the group $\Pi$ admits an embedding $\Pi \hookrightarrow (F_2)^{\times 4}$ via an inclusion of the general fiber $U_t \hookrightarrow \mathcal{U}=(E\setminus \{0\})^{\times 4}$ and the projection of $\Pi$ on each factor $F_2$ is surjective.

Thanks to Proposition \ref{local factor prop} and Corollary \ref{triple cover lemma}, we are able describe the kernel of $\Pi \to \Gamma$ pretty explicitly. Let $F_i=\langle a_i, b_i \rangle, \ i-1, \ldots, 4,$ be the $i$-th factor in $(F_2)^{\times 4}=\pi_1(\mathcal{U})$. Then a loop around the $i$-th boundary component in $\mathcal{U}=\prod_{i=1}^{4}(E_i \setminus \{0\})$ represents a conjugate of $c_i:=[a_i, b_i]$.

Proposition \ref{local factor prop} gives an explicit description of the kernel of the map 
\[
\pi_1(V_0(p)) \to \pi_1(\widehat{V}(p)).
\]
At the point $p$ the three curves $B_0, \ B_1, \ B_2$ meet pairwise transversally and the fundamental group of $\pi_1(V_0(p))$ is generated by loops $\mu_0, \ \mu_1, \ \mu_2$ around them with $\epsilon=\mu_0\mu_1\mu_2$ as a central element. The kernel $\pi_1(V_0(p)) \to \pi_1(\widehat{V}_0(p))$ is normally generated by the element $\rho_P:=\mu_1\epsilon=\mu_1\mu_0\mu_1\mu_2$. The images of the loops $\mu_0, \ \mu_1, \ \mu_2$ in $\Pi$ under the map $\pi_1(U_0) \to \Pi$ are some elements $(m_1,1,1,1), \  (1,m_2,1,1)$ and $(1,1,1, m_4)$, where  $m_i$ are $F_i$-conjugates of $c_i=[a_i,b_i], \ i=1, 2, \ 4.$ respectively. The image of the element $\rho_P$ equals $r_P=(m^2_1,m_2,1,m_4)$ (recall that $c_i$ pairwise commute in $\Pi$, and so do their conjugates).

Similarly, the local fundamental group of $U_0$ around $q$ is generated by loops $\nu_1,  \ \nu_2,  \ \nu_{\tau}$ around the divisors $B_1, \ B_2, \ B_{\tau}$ with a central element $\xi=\nu_1\nu_2\nu_{\tau}$ and the kernel generated by $\rho_Q:=\nu_2\xi=\nu_2\nu_1\nu_2\nu_{\tau}$. In $\Pi$ this element projects to $r_Q:=(n_1,n_2^2, n_3, 1)$. 

Summarizing all this, we can get a relatively explicit description of the group $\Gamma$.

\begin{prop}\label{pi descrip prop}
Let $F_i, \ i=1, \ldots, 4$ be copies of the free group on two generators, $F_i=\langle a_i, b_i \rangle$. Denote $c_i:=[a_i, b_i]$ and $H=H_i:=F_i/[F_i, F_i]$. Let $\Pi$ be the kernel of the map
\[
\prod_i F_i \to \bigoplus H_i \xrightarrow{\Phi} H^{\oplus 2},
\]
where the first arrow is the abelianisation and $\Phi$ is given by an integral matrix
\[
\Phi=\begin{pmatrix} \Id_H& \Id_H& -\Id_H& -\Id_H \\\Id_H& 0& \Id_H& -2\cdot \Id_H \end{pmatrix}
\]

Then
\[
\Gamma \simeq \Pi / \langle \langle r_P, r_Q \rangle \rangle,
\]
where $r_P=(m^2_1,m_2,1,m_4)$ and $r_Q=(n_1,n^2_2,n_3,1)$ are elements of $\prod_i F_i$, where  $m_i$ and $n_i$ are certain $F_i$-conjugate of $c_i$.
\end{prop}

Remark that  $r_P$ and $r_Q$ are automatically contained in $\Pi$ because both go to zero in the abelianisation. 

Recall that the group $\Pi$ contains a normal subgroup $K \simeq \prod([F_i, F_i])$ which is the kernel of the map $\Pi \to \Pi^{\operatorname{ab}}$.

\begin{prop}\label{metabelian prop}
Let $\pi \colon \Pi \to \Gamma$ be the projection. The image $M:=\pi(K)$ is abelian in $\Gamma$. It is normal in it and $L:=\Gamma/M$ is isomorphic to $\Z^4$. In particular, $\Gamma$ is metabelian.
\end{prop}
\begin{proof}
The group $K$ is generated by mutually commuting subgroups $DF_i=[F_i, F_i]$. It is enough to show that for each $i$ the group $\pi(DF_i)$ is abelian.

Note that for each $i$ the projection $\Pi \to F_i$ is surjective because $\Pi$ is the full preimage of $L=\ker \Phi$ in $\prod_i F_i^{\ab}$ and $p_i \colon L \to H_i=F_i^{\ab}$ is onto for every $i$. 

Each group $DF_i$ is generated by $F_i$-conjugates of $c_i$. The argument is  the same for every $i$ up to  small adjustments depending on whether $i\in\{2,4\}$ or $i \in \{1, 3\}$.

Assume the first case. We claim that $\pi(m_i)$ is a central element in $\pi(DF_i)$. Indeed,for every $x \in DF_i$ one has
\[
1=[\pi(r_P), \pi(x)]=[\pi(m_i), \pi(x)].
\]
At the same time, $m_i$ is an $F_i$-conjugate of $c_i$, so $DF_i$ is generated by conjugates of $c_i$, or, equivalently, by elements of the form $m'_i=ym_iy^{-1}, \ y \in F_i$. As we point out above, the projection $\Pi \to F_i$ is surjective, so  such $y$ admits a lift $\widetilde{y} \in \Pi$ which allows  us to realise $m'_i$ as a $\Pi$-conjugate of $m_i$. By a similar argument
\[
1=[\pi(\widetilde{y}r_P\widetilde{y}^{-1}),\pi(x)]=[\pi(m'_i), \pi(x)]
\]
we conclude that $\pi(m'_i)$ is central too. In other words, $\pi(DF_i)$ is generated by central elements, hence commutative.

In the case where $i \in \{1,3\}$ the argument is the same with using $r_Q$ instead of $r_P$ and $n_i$ instead of $m_i$. 

This proves that $M$ is abelian.

Since $K$ is normal in $\Pi$ and $\pi$ is surjective, $M$ is normal in $\Gamma$. Finally, we obtain a commutative diagram
\[
\xymatrix{
1 \ar[r]& K \ar[r] \ar[d]^{\pi} & \Pi \ar[d]^{\pi} \ar[r]& \ker \Phi\simeq \Z^4 \ar[d]^{\sigma} \ar[r] &1\\
1 \ar[r]& M \ar[r] & \Gamma \ar[r]& \Gamma/M \ar[r] &1
}
\]
Since $r_P, r_Q \in K$ the kernel of $\pi$ is in $K$ and a simple diagram chasing shows  that $\sigma$ is an isomorphism.

\end{proof}

Combining with Proposition \ref{metab kill prop}, we obtain the following corollary.

\begin{cor}
$\Gamma$ is virtually nilpotent.

\end{cor}

\subsection{$\pi_1(Y)$ is rationally 3-step nilpotent}\label{subsec last}

In this subsection we combine the description of $\Gamma$ obtained in subsection \ref{metabelian subsec} with the criteria of subsection \ref{group prep}. Recall that $\Gamma$ fits into an exact sequence
\[
1\to M\to\Gamma\to L\to 1,\qquad M=\pi(K),\quad L\simeq\mathbb Z^4,
\]
and that by Proposition \ref{pi descrip prop} one has $\Gamma\simeq\Pi/\langle\langle r_P,r_Q\rangle\rangle$, where
\[
r_P=(m_1,m_2^2,1,m_4),\qquad r_Q=(n_1^2,n_2,n_3,1).
\]
Our task is to describe $V:=M\otimes\mathbb Q$ as a $\mathbb Q[L]$-module. Recall that the lattice $L$ is given as the kernel of the matrix $\Phi = \overline{f}_{*}$ (see equation (\ref{matrix})). Solving $\Phi(x_1, \ldots ,x_4)=0$ one gets
\begin{equation}\label{L formula}
L:=\{(2h'-h, \  2h-h', \  h,  \ h' )\  | h, h' \in H \} \simeq H \oplus H
\end{equation}
(cf. proof of  Proposition \ref{exact app prop}). Denote by $p_i \colon L \to H_i$ the $i$-th projection of $L$ onto $H_i=\pi_1(E_i)$, i.e. $p_1(h, h')=2h'-h, \ p_2(h,h')=2h'-h$ etc.

The group $\Pi$ is an extension of $L$ by $K:=\prod_i (DF_i)$. This gives an action of the ring $\Z[L]$ on $K^{\ab}$. Let us denote $S:=\Z[H]$. We also sometimes write $S_i$ for $S$ viewed as a $\Z[L]$-module where the action is induced by $p_i \colon L \to H$.

\begin{prop}
Let $[c_i]$ be the class of $c_i \in DF_i$ in the abelianisation $DF_i^{\ab}$.  Then there is an isomorphism of $\Z[L]$-modules
\[
K^{\ab} \simeq \bigoplus_{i=1}^4 S_i\cdot [c_i].
\]
\end{prop}
\begin{proof}
To see how an element $\ell \in L$ acts on $K^{\ab}$ choose a lift $\widetilde{\ell}$ to $\Pi$ and let it act by conjugation on $\Pi/[K,K]$. This action restricts from the conjugation action of $\prod (F_i)^{\ab}=\prod H_i$ on $\prod (DF_i)^{\ab}$. Each $H_i$ acts on $(DF_j)^{\ab}$ trivially for $j \neq i$ because different factors in $\prod_i F_i$ commute. It is enough to verify that for a free group $F$ on $2$ generators $(DF)^{\ab}$ is a free $\Z[F^{\ab}]$-module generated by the class $\overline{c}$ of $c=[a, b]$. To see this, represent $F$ as the fundamental group of a punctered elliptic curve $E \simeq \{0\}$. Then $(DF)^{\ab}$ identifes with the first homology group of the maximal abelian cover of $E\setminus \{0\}$ that is homeomorphic to $\R^2 \setminus \Z^2$. This homology group is freely generated by loops around the punctures and $\pi_1(\R^2/\Z^2)$ acts on it simply transitively.
\end{proof}

Recall that the kernel of the epimorphism $\pi \colon \Pi \to \Gamma$ is normally generated by elements $r_P=(m_1^2, m_2, 1, m_4)$ and $r_Q=(n_1, n_2^2, n_3, 1)$, where $m_i$ and $n_i$ are conjugates of $c_i$ in $F_i$. Write $m_i=y_ic_iy_i^{-1}$ and $n_i=z_ic_iz_i^{-1}$ and denote by $\overline{y_i}$ and $\overline{z_i}$ their classes in $H_i=F_i^{\ab}$. We also denote by $[c_i], \ [r_P], \ [r_Q]$ the images of $c_i, \ r_P, r_Q$ in $(DF_i)^{\ab}$. 

\begin{prop}\label{module desc prop}
There is an isomorphism of $\mathbb Z[L]$-modules
\[
M\simeq\Big(\bigoplus_{i=1}^4S_i[c_i]\Big)\Big/\,\mathbb Z[L]\cdot\langle [r_P], [r_Q]\rangle,
\]
where
\begin{equation}\label{rel}
[ r_P]=2\bar y_1[c_1]+\bar y_2[c_2]+\bar y_4[c_4],\qquad
[r_Q]=\bar z_1[c_1]+2\bar z_2[c_2]+\bar z_3[c_3].
\end{equation}
\end{prop}

\begin{proof}
This is the abelianisation of the presentation from Proposition \ref{pi descrip prop}. Since $M$ is abelian (Proposition \ref{metabelian prop}), $[K,K]$ lies in the normal closure of $r_P, r_Q$, which implies the claim.
\end{proof}

From now on we work rationally. We will use the presentation of $M$ from Proposition \ref{module desc prop} to study the $\Q[L]$-module $V:=M \o \Q$ and apply the criterion from Lemma \ref{rational criterion}.

Let $\mathcal{S}:=\bigoplus_i S_i [c_i] \o \Q$ and $\mathcal{R}:=\Q[L]\langle[r_P], [r_Q] \rangle$, so that $V=\mathcal{S}/\mathcal{R}$.

According to the strategy of Lemma \ref{rational criterion}, we need to understand the filtration of $V$ given by the powers of the augmentation ideal $I \subset \Q[L]$. 

\begin{prop}\label{ann prop}
$V$ is annihilated by some power of $I$.
\end{prop}
\begin{proof}
By Proposition \ref{triv supp prop} it is enough to show that $V_{\C}$ is supported on a trivial character, i.e. $V^{ss}$ is trivial. Let $T_i:=\im[p_i^* \colon \mathbb{T}(H) \to \mathbb{T}(L)]$. Since $V$ is a quotient of $\mathcal{S}$, one has
\[
\operatorname{Supp(V_{\C})}\subseteq \operatorname{Supp}(\mathcal{S_{\C}})=\bigcup_i T_i.
\]

Assume that $\operatorname{Supp}(V_{\C})$ contains a non-trivial character $\chi$; we derive a contradiction by looking at  which of the subtori $T_i$ contain this character. Fisrt, $\chi$  cannot be contained in exactly one subtorus, as we now explain.

Denote $S_{i, \chi}:=(S_i \o \C)\o_{\C[L]} \C_{\chi}$. If a character $\chi$ is contained in $T_i$ but not in $T_j$ then $\dim S_{i, \chi}=1$ and $S_{j, \chi}=0$. Indeed, since $S_i$ is a free cyclic $H$-module and $p_i \colon L \to H$ is surjective, $S_{i, \chi}=\C\chi/(\ker p_i\cdot \C_{\chi})$ and the question reduces to whether there exists $\ell \in L$ such that $p_i(\ell) = 0$ and $\chi(\ell) \neq 1$. If such $\ell$ does not exist, i.e. $\chi$ factors through $p_i$, then $\ker p_i \cdot \C_{\chi}$ is trivial. Otherwise, $\ell-1$ acts non-trivially on $\C_{\chi}$ and $S_{\chi, i}$ is trivial.

Suppose that $\chi \in \operatorname{Supp}(V)$ is contained in exactly one subtorus $T_i$. Then the module $(\mathcal{S}_{\C}/\mathcal{R}_{\C})\o_{\C[L]} \C_{\chi}$ is a quotient of the line $S_{i, \chi}$ spanned by $[c_i]$. By (\ref{rel}), the image of $[r_P]$ (if $i \in \{2,4\}$) or $[r_Q]$($i \in \{1,3\}$) is a non-zero multiple of $[c_i]$. Therefore $(\mathcal{S}_{\C}/\mathcal{R}_{\C}) \o_{\C[L]} \C_{\chi}=0$.

Using the splitting (\ref{L formula}), decompose $\mathbb{T}(L)=\mathbb{T}(H) \times \mathbb{T}(H), \ \chi=(\chi_1, \chi_2)$. With respect to this splitting, $T_1=\{\chi_1^{2}=\chi_2^{-1}\}, \ T_2=\{\chi_2^2=\chi_1^{-1}\}, \ T_3=\{\chi_1=1\}$ and $T_4=\{\chi_2=1\}$.  Observe that the pairwise intersections $T_i \cap T_j$ are torsion for each $i$ and $j$ and $T_3 \cap T_4$ is trivial.

Suppose that $\chi \in T_1 \cap T_3$ is a non-trivial character in the support of $V_{\C}$.The fiber of $\mathcal{S}_{\C}$ at $\chi$ is spanned by $[c_1], [c_3]$ and the relations are given by the vectors $\chi([r_P])=(2\chi(\bar y_1), 0)$ and $\chi([r_Q])=(\chi(\bar z_1), \chi(\bar z_3))$ . The matrix $\begin{pmatrix} 2\chi(\bar y_1) & \chi(\bar z_1) \\ 0 & \chi(\bar z_3)\end{pmatrix}$ is non-degenerate, which means that $\dim \langle \chi([r_P]), \chi([r_Q]) \rangle=2$  and the quotient of the two-dimensional space $\langle[c_1], [c_2] \rangle$ by this space is zero.

The same argument works if $\chi$ is in either $T_1 \cap T_4, \  T_2 \cap T_3$ or $T_2 \cap T_4$: in all three cases the matrix of relations is either upper or lower triangular. Finally, for $\chi \in T_1 \cap T_2$ the matrix of relations has form $\begin{pmatrix} 2 \chi(\bar y_1) & \chi(\bar z_1) \\ \chi(\bar y_2) & 2\chi(\bar z_2)  \end{pmatrix}$. The tori $T_1$ and $T_2$ intersect in $3$-torsion, therefore the elements of this matrix are roots of unity and its determinant  $4\chi(\bar y_1)\chi(\bar z_2)-\chi(\bar y_2)\chi(\bar z_1)$ has absolute value at least $3$. Hence, this matrix is non-degenerate is well.

\end{proof}

Now we are able to apply Lemma \ref{rational criterion} and compute the dimensions of rational lower central series of $\Gamma$ explicitly.

\begin{prop}\label{i2 prop}\label{i3 prop}
One has
\[0=I^2V \subset IV \subset V\]
with $\dim V/IV=\dim IV/I^2V=2$.
\end{prop}
\begin{proof}
We have $V/I^kV=\mathcal{S}/(I^k\mathcal{S}+\mathcal{R})$.
\\

\textbf{Step 1. Calculating $ V/IV$}. Observe that $I^kS_i=J^kS_i$, where $J \subseteq S$ is the augmentation ideal. Indeed, the action of $L$ on $S$ comes from the ring homomorphism $p_i \colon \Z[L] \to S_i$ and $p_i(I) \subseteq J$. At the same time, $p_i$ is surjective, as one can easily see from the formula (\ref{L formula}), thus every element $x \in J$ has a preimage in the augmentation ideal.

For $k=1$ the space $V/IV$ is the quotient of  $\langle[c_1], [c_2], [c_3], [c_4] \rangle_{\Q}$ by the subspace generated by the vectors $(2,1,0,1)$ and $(1, 2, 1,0)$. These two vectors are linearly independent and $\dim V/IV=2$.
\\

\textbf{Step 2. Calculating $V/I^2V$.} Let us now understand $V/I^2V$. By the discussion above, $\mathcal{S}/I^2\mathcal{S} = \bigoplus_{i} (S/J^2S\cdot[c_i])$. The ring $S =\Z[H]$ is isomorphic to the ring of Laurent polynomials in two variables, $S \simeq \Z[u^{\pm 1}, v^{\pm1}]$.

The augmentation ideal $J$ is the maximal ideal of the point $(1,1) \in \operatorname{Spec}(S)=\mathbb{G}_m \times \mathbb{G}_m$. Thus, $\dim S/IS=1$ and $\dim (IS/I^2S)=\dim T_{(1,1)}\operatorname{Spec}(S)= 2$. One concludes that $\dim S/I^2S=3$ and $\dim\mathcal{S}/I^2\mathcal{S}$ is four times bigger, that is, $12$.

The submodule $\mathcal{R} \subseteq \mathcal{S}$ is generated by $[r_P]$ and $[r_Q]$ over $L$. Modulo $I^2\mathcal{S}$ it is linearly generated by $[r_P], [r_Q]$ and elements of the form $(\ell-1) \cdot[r_P], \ (\ell-1) \cdot[r_Q]$. Observe that 
\[
(p_i(\ell)-1)\bar y_i=(p_i(\ell)-1)(1+(\bar{y}_i-1))\equiv (p_i(\ell)-1) \ \operatorname{mod}J^2S_i.
\]
This implies that
\begin{align*}
(\ell-1)\cdot[r_P]
&=2(p_1(\ell)-1)\,\bar y_1[c_1]+2\,(p_2(\ell)-1)\,\bar y_2[c_2]+(p_4(\ell)-1)\,\bar y_4[c_4]\\
&\equiv 2p_1(\ell)[c_1]+p_2(\ell)[c_2]+p_4(\ell)[c_4]\pmod{I^2\mathcal S}.
\end{align*}
Similarly, 
\[
(\ell-1)\cdot [r_Q]\equiv 2p_1(\ell)[c_1]+p_2(\ell)c_2+p_3(\ell)[c_3] \pmod{ I^2\mathcal{S}}.
\]
Varying $\ell$, we get two vector spaces, $A_P:=I\cdot[r_P]/I^2\mathcal{S}$ and $A_Q:=I \cdot[r_Q]/I^2\mathcal{S}$. Since $p_1 \oplus p_2 \colon L \to H \oplus H$ is injective, they both have dimension at least $4$. Moreover, from formula (\ref{L formula}) one has $p_3=\frac{p_1+2p_2}{3}$ and $p_4=\frac{2p_1+p_2}{3}$, thus $A_P$ and $A_Q$ have dimension exactly $4$.

One has $\mathcal{R}/(\mathcal{R} \cap I^2\mathcal{S})=\Q[r_Q] \oplus \Q[r_P]\oplus(A_P+A_Q)$, so in order to compute its dimension we need to calculate $\dim(A_P\cap A_Q)=\dim A_P+\dim A_Q -\dim(A_P+A_Q)$.  Assume $a \in A_P \cap A_{Q}$. Then 
\[
a=(2p_1(\ell),p_2(\ell), 0, p_4(\ell))=(p_1(\ell'), 2p_2(\ell'), p_3(\ell'), 0).
\]
Again, using description (\ref{L formula}) one checks that this is only possible when $\ell=(h,0)$ and $\ell'=(0, -h)$ for some $h \in H\o\Q$; then $a=(-2h,h, 0,0)$. In other words, 
\[
\dim_{\Q}(A_P\cap A_Q)=\dim_{\Q} H_{\Q}=2.
\]
 We conclude that $\dim (\mathcal{R}/(I^2\mathcal{S}\cap \mathcal{R}))=2+4+4-2=8$ and $\dim V/I^2V=12-8=4$. Therefore $\dim IV/I^2V=\dim V/I^2V-\dim V/IV=2.$
\\

\textbf{Step 3. Calculating $V/I^3V$}. Let us show that $I^2\mathcal{S} \subseteq \mathcal{R}+I^3\mathcal{S}$. This would imply that $I^2V=I^3V$, and since $V$ is annihilated by some power of $I$ (Proposition \ref{ann prop}) this is enough to conclude that $I^3V=0$.

Set $\widetilde{[r_P]}=2[c_1]+[c_2]+[c_4]$ and $\widetilde{[r_Q]}=[c_1]+2[c_2]+[c_3]$. Then $[r_P]-\widetilde{[r_P]}$ and $[r_Q]-\widetilde{[r_Q]}$ are contained in $I\mathcal{S}$, which means that $s \cdot [r_P] \equiv s\cdot \widetilde{[r_P]} \pmod {I^3\mathcal{S}}$ and $s\cdot [r_Q] \equiv s\cdot \widetilde{[r_Q]} \pmod{ I^3\mathcal{S}}$ for any $s \in I^2 \subseteq \Q[L]$.  As we show below, this allows us to find, for every $i$, a multiple of either $r_P$ or $r_Q$ that becomes proportional to $[c_i]$ modulo $I^3\mathcal{S}$

Choose $\ell, \ \ell' \in L$ such that 
$p_2(\ell)=p_3(\ell')=0$. Then 
$(\ell-1)(\ell'-1) r_Q$ is proportional to $[c_1]$, namely 
\[
(\ell-1)(\ell'-1) [r_Q]\equiv(\ell-1)(\ell'-1)\widetilde{[r_Q]}=(p_1(\ell)-1)(p_1(\ell')-1)[c_1] \pmod{I^3\mathcal{S}}.
\]
We would like to know the answer to the converse problem: when for given $a, b \in H$ can we realize $(a-1)(b-1)[c_1]$ as $(\ell-1)(\ell'-1)[r_Q]$ for some $\ell, \ \ell' \in L$. A difficulty here is that $p_1(\ker p_2)=3H$ and $p_1(\ker p_3)=2H$, so in  such elements do not generate the module $I^2S_1$ integrally. Nevertheless, they do generate it rationally modulo $I^3S_i$. The reason behind this is that for a natural number $n$ and an element $a \in H$ one has 
\[
n\cdot a-1=a^{n}-1=(a-1)(1+a+\ldots+a^{n-1}) \equiv n(a-1) \pmod{J^2}
\]
(here, to avoid any confusion, $n \cdot a$ denotes multiplication by $n$ in the group, not in the group ring). Thus, $(2a-1)(3b-1)\cdot[c_1] \equiv 6(a-1)(b-1)[c_1] \pmod{J^3S_1}$, so such elements generate $J^2/J^3$ over $\Q$. Arguing in the same spirit (substituting $r_P$ where necessary), one shows that $I^2S_i[c_i] \subseteq I^2\mathcal{R} +I^3S_i$ for every $i$. As we explain above, this implies the statement.

\end{proof}

\begin{thrm}\label{final thrm}
The group $\Gamma=\pi_1(Y)$ is virtually nilpotent of rational nilpotency class $3$.The same holds for every finite index subgroup of $\Gamma$. In particular, $\Gamma$ is not virtually $2$-step nilpotent.
\end{thrm}

\begin{proof}
Let $W:=\gamma_1\Gamma \o \Q$. By Lemma \ref{rational criterion} $\gamma_k \Gamma \o \Q=I^kW$. It is therefore enough to show that $W=V$, i.e.\ that the map $\eta\colon\Lambda^2L\to V/IV$ is rationally surjective.

Let $\omega$ be the intersection form on $H$. For $x,x'\in F_i$ one has $[x,x']\equiv c_i^{\,\omega(\bar x,\bar x')}$ modulo $[F_i,DF_i]$, and the image of $[F_i,DF_i]$ in $DF_i^{ab}$ is $I_H\,DF_i^{ab}$. Hence, for $\ell,\ell'\in L$ with lifts $\tilde\ell,\tilde\ell'\in\Pi$, the class of $[\tilde\ell,\tilde\ell']\in K$ modulo $I\,K^{ab}$ equals
\[
\sum_{i=1}^4\omega\big(\mathrm{pr}_i\ell,\mathrm{pr}_i\ell'\big)\,[c_i] .
\]
Choose $h,h'\in H$ with $\omega(h,h')=1$. The pairs $\big((h,0),(h',0)\big)$, $\big((0,h),(0,h')\big)$ and $\big((h,0),(0,h')\big)$ give the vectors
\[
(1,4,1,0),\qquad (4,1,0,1),\qquad (-2,-2,0,0)
\]
in $\mathbb Q^4=\langle [c_1],\dots,[c_4]\rangle$. Together with $\langle(2,1,0,1), (1,2,1,0)\rangle=\im[\mathcal{R} \to \mathcal{S}]$ they span $\mathbb Q^4$. In the notations of Lemma \ref{rational criterion} we conclude that $\eta$ is onto $V/IV$, and $W=V$. The claim then follows from Lemma \ref{rational criterion} and Proposition \ref{i3 prop}.

\end{proof}
\bibliography{references.bib}{}

\begin{thebibliography}{CDHP26}

\bibitem[AAC25]{AC}
Rodolfo Aguilar~Aguilar and Fr{\'e}d{\'e}ric Campana.
\newblock The nilpotent quotients of normal quasi-projective varieties with proper quasi-{A}lbanese map.
\newblock {\em Pure Appl. Math. Q.}, 21(3):911--929, 2025.

\bibitem[AN99]{AN}
Donu Arapura and Madhav Nori.
\newblock Solvable fundamental groups of algebraic varieties and {K}{\"a}hler manifolds.
\newblock {\em Compositio Math.}, 116:173--188, 1999.

\bibitem[Cam95]{Campana95}
Fr{\'e}d{\'e}ric Campana.
\newblock Remarques sur les groupes de {K}{\"a}hler nilpotents.
\newblock {\em Ann. Sci. {\'E}cole Norm. Sup. (4)}, 28(3):307--316, 1995.

\bibitem[CDHP26]{CDHP}
Junyan Cao, Ya~Deng, Christopher~D. Hacon, and Mihai P{\u{a}}un.
\newblock Two-step nilpotent monodromy of local systems on special varieties, 2026.
\newblock Preprint, arXiv:2603.14539.

\bibitem[CT95]{CT}
James~A. Carlson and Domingo Toledo.
\newblock Quadratic presentations and nilpotent {K}{\"a}hler groups.
\newblock {\em J. Geom. Anal.}, 5(3):359--377, 1995.

\bibitem[Hat75]{Hat}
Akio Hattori.
\newblock Topology of $\mathbb{C}^n$ minus a finite number of affine hyperplanes in general position.
\newblock {\em J. Fac. Sci. Univ. Tokyo Sect. IA Math.}, 22(2):205--219, 1975.

\bibitem[Mat12]{Mather}
John Mather.
\newblock Notes on topological stability.
\newblock {\em Bull. Amer. Math. Soc. (N.S.)}, 49(4):475--506, 2012.
\newblock Reprint of lecture notes, Harvard University, 1970.

\bibitem[Nor83]{Nori}
Madhav~V. Nori.
\newblock Zariski's conjecture and related problems.
\newblock {\em Ann. Sci. {\'E}cole Norm. Sup. (4)}, 16(2):305--344, 1983.

\bibitem[Py25]{Py}
Pierre Py.
\newblock {\em Lectures on {K}{\"a}hler Groups}.
\newblock Princeton Mathematical Series. Princeton University Press, Princeton, NJ, 2025.

\bibitem[Rog24]{rogbns}
Vasily Rogov.
\newblock The {B}ieri--{N}eumann--{S}trebel sets of quasi-projective groups, 2024.
\newblock Preprint, arXiv:2408.06250 [math.AG].

\bibitem[Rog25]{Rogomin}
Vasily Rogov.
\newblock o-minimal geometry of higher {A}lbanese manifolds, 2025.
\newblock Preprint, arXiv:2505.07632 [math.AG].

\bibitem[Rog26]{Rogcd}
Vasily Rogov.
\newblock Topology of higher {A}lbanese maps and aspherical varieties with nilpotent fundamental group, 2026.
\newblock Preprint, arXiv:2607.19013.

\bibitem[Shi25]{Shim2}
Taito Shimoji.
\newblock On the torsion-free nilpotent fundamental groups of smooth quasi-projective varieties of rank up to seven, 2025.
\newblock Preprint, arXiv:2510.09026 [math.AG].

\bibitem[Shi26]{Shim1}
Taito Shimoji.
\newblock Gradings on nilpotent lie algebras associated with the nilpotent fundamental groups of smooth complex algebraic varieties.
\newblock {\em Journal of Pure and Applied Algebra}, 230(1):108154, January 2026.

\end{thebibliography}
\bibliographystyle{alpha}
\end{document}